\documentclass[11pt]{amsart}
\usepackage{mathptmx}

\usepackage{amsmath}
\usepackage{amsthm, amssymb}
\usepackage{mathbbol}
\usepackage[ansinew]{inputenc}
\usepackage{graphicx}
\usepackage{color}

\newtheorem{thm}{Theorem}
\newtheorem{prop}[thm]{Proposition}
\newtheorem{lemma}[thm]{Lemma}
\theoremstyle{definition}
\newtheorem{defi}[thm]{Definition}

\theoremstyle{remark}
\newtheorem{ex}[thm]{Example}
\newtheorem{rem}[thm]{Remark}

\newcommand{\alg}[1]{\mathfrak{{#1}}}
\newcommand{\tr}{\text{tr}}
\newcommand{\cn}[1]{ {\mathbb{C}^{{#1}}} }
\newcommand{\rn}[1]{ {\mathbb{R}^{{#1}}} }
\newcommand{\co}[2]{ {\left[{#1},{#2}\right]} } 
\newcommand{\aco}[2]{ {\left\{{#1},{#2}\right\}} }
\newcommand{\eref}[1]{(\ref{#1})} 

\newcommand{\pderi}[2]{ { \frac{\partial {#1} }{\partial {#2} } } }
\newcommand{\ad}{{\text{ad}}}

\begin{document}
\title{A Counterexample to the Quantizability of Modules}
\author{Thomas Willwacher}
\address{Department of Mathematics, ETH Zurich}
\email{thomas.willwacher@math.ethz.ch}
\thanks{The author was partially supported by the Swiss National Science Foundation (grant 200020-105450)}
\date{}

\begin{abstract}
Let $\pi$ be a Poisson structure on $\rn{n}$ vanishing at 0. It leads to a Kontsevich type star product $\star_\pi$ on $C^\infty(\rn{n})[[\epsilon]]$. We show that 
\begin{enumerate}
\item The evaluation map at 0 
\[
ev_0: C^\infty(\rn{n})\rightarrow \mathbb{C}
\]
can in general not be quantized to a character of $(C^\infty(\rn{n})[[\epsilon]],\star_\pi)$. 
\item A given Poisson structure $\pi$ vanishing at zero can in general not be extended to a formal Poisson structure $\pi_\epsilon$ also vanishing at zero, such that $ev_0$ can be quantized to a character of $(C^\infty(\rn{n})[[\epsilon]],\star_{\pi_\epsilon})$. 
\end{enumerate}
We do not know whether the second claim remains true if one allows the higher order terms in $\epsilon$ to attain nonzero values at zero. 
\end{abstract}
\maketitle

\section*{How to read this paper in 2 minutes}
The busy reader can take the following shortcut:
\begin{enumerate}
\item Read Theorem \ref{thm:main} on page \pageref{thm:main} for the main result.
\item Read Definition \ref{def:evmquant} if its statement is not clear.
\item Look at eqns. \eref{equ:cexdef} and the preceding enumeration for the definition of the counterexample.
\end{enumerate}

\section{Introduction}

Let $M$ be a smooth $n$-dimensional manifold equipped with a Poisson structure $\pi$, making $C^\infty(M)$ a Poisson algebra with bracket $\aco{\cdot}{\cdot}$. In this paper, we will exclusively deal with the case $M=\rn{n}$. Kontsevich \cite{kontsevich} has shown that one can always quantize this algebra, i.e., find an associative product $\star_\pi$ on $C^\infty(M)[[\epsilon]]$ such that for all $f,g\in C^\infty(M)$
\[
f\star_\pi g = fg + \frac{\epsilon}{2} \aco{f}{g} + O(\epsilon^2).
\]  
Furthermore, he showed that the set of such star products is, up to equivalence, in one to one correspondence with the set of formal Poisson structures on $M$ extending $\pi$. 
\begin{defi}
A \emph{formal Poisson structure} $\pi_\epsilon$ is a formal bivector field 
$$
\pi_\epsilon \in \Gamma (\Lambda^2TM)[[\epsilon]]
$$
satisfying the Jacobi identity
\begin{equation}
\label{equ:piepsjacobi}
\co{\pi_\epsilon}{\pi_\epsilon}=0
\end{equation}
where $\co{\cdot}{\cdot}$ is the Schouten-Nijenhuis bracket. We say that $\pi_\epsilon$ \emph{extends} the Poisson structure $\pi\in \Gamma(\Lambda^2TM)$ if its $\epsilon^0$-component is $\pi$.
\end{defi}  

Let now $m\in M$ be a point and consider the evaluation map
$$
ev_m:C^\infty(M)\ni f\mapsto f(m) \in \cn{}.
$$
It makes $\cn{}$ into a $C^\infty(M)$-module, i.e., for all $f,g\in C^\infty(M)$
$$
ev_m(fg)=ev_m(f)ev_m(g).
$$
The question treated in this paper is the following:

\vspace{5mm}
\textbf{Main Question:} Can one quantize the evaluation map $ev_m$?
\vspace{5mm}

By this we mean the following:
\begin{defi}
\label{def:evmquant}
Let $\star_{\pi_\epsilon}$ be the Kontsevich star product associated to the formal Poisson structure $\pi_\epsilon$ on $M$ and let $m\in M$ be an arbitrary point.
A formal map 
$$
\rho: C^\infty(\rn{n}) \rightarrow \mathbb{C} [[ \epsilon ]]
$$
will be called \emph{quantization of $ev_m$} if the following holds.
\begin{itemize}
\item It has the form
\[
\rho(f) = f(m) + \epsilon \rho_1(f) + \epsilon^2 \rho_2(f) + \cdots
\]
where the $\rho_k$ are differential operators evaluated at $m$. Concretely, this means in local coordinates that
\[
\rho_k(f) = \sum_I c_I \pderi{f}{x^I}(m)
\]
where the sum is over multiindices and the $c_I$ are constants, vanishing except for finitely many $I$.
\item For all $f,g\in C^\infty(M)$
\begin{equation}
\label{equ:mod}
\rho(f\star_{\pi_\epsilon} g) = \rho(f) \rho(g)
\end{equation}
\end{itemize}
\end{defi}

\begin{lemma}
Let $\pi_\epsilon=\pi+O(\epsilon)$ and $m\in M$. If a quantization of $ev_m$ exists, then $\pi(m)=0$, i.e., $\pi$ vanishes at $m$.
\end{lemma}
\begin{proof}
The $\epsilon^1$-component of the equation
$$
\rho(\co{f}{g}_\star) = \co{\rho(f)}{\rho(g)} = 0
$$ 
reads
$$
\aco{f}{g}(m) = 0.
$$
Hence $\pi(m) = 0$.

Remark: A similar calculation for a higher dimensional submanifold also yields the higher dimensional coisotropy condition.
\end{proof}

From now on we will assume that $\pi(m) = 0$, or, equivalently, that $\{m\}\subset M$ is coisotropic. For details on coisotropic submanifolds see \cite{cattaneo-2004-69}.

The above main question has been answered positively by Cattaneo and Felder in \cite{cattaneo-2004-69}, \cite{cattaneo-2005}, provided $\pi$ satisfies certain conditions. Adapted to our context, they proved the following theorem.

\begin{thm}
\label{thm:cf}
For any formal Poisson structure $\pi_\epsilon$ on $M=\rn{n}$ such that
\[
\pi_\epsilon(m) = 0 + O(\epsilon)
\] 
for some $m\in M$ there exists a linear map
\[
\tilde{\rho}: C^\infty(M) \rightarrow \mathbb{C} [[ \epsilon ]]
\]
\[
\tilde{\rho}(f) = f(m) + \epsilon^2 \tilde{\rho}_2(f) + \cdots \footnote{Note the absence of the $\epsilon^1$ term.}
\]
where the $\tilde{\rho}_k$ are differential operators evaluated at $m$, that satisfies
\begin{equation}
\label{equ:modano}
\tilde{\rho}(f\star_{\pi_\epsilon} g) = \tilde{\rho}(f) \tilde{\rho}(g) + A(f,g)
\end{equation}
for all $f,g\in C^\infty(M)$. Here $A=O(\epsilon^2)$ (the ``anomaly'') is a bidifferential operator.
\end{thm}

Hence, if the ``anomalous'' term on the r.h.s. of \eref{equ:modano} vanishes, one sees that $\tilde{\rho}$ becomes a quantization of $ev_m$. The precise form of $A(f,g)$ is recalled in section \ref{sec:cfano}.

When $\pi_\epsilon(m) = 0$ to all orders in $\epsilon$, the anomaly is actually at least of order $\epsilon^3$. Furthermore, we will later provide an example for which the $\epsilon^3$-term does not vanish.

A theorem similar to Theorem \ref{thm:cf} above holds in the case of higher dimensional coisotropic submanifolds. There, anomaly terms will also occur in general. It is still an open question whether the vanishing of these terms is merely a removable technical condition or a fundamental obstruction to quantizability. Our paper gives a partial answer to this question in the simplest possible case.

\subsection{Quantization of Modules}
\label{sec:cafmod}
In this section we review the construction of Cattaneo and Felder \cite{cattaneo-2004-69} leading to Theorem \ref{thm:cf}.
We throughout assume familiarity with the construction of Kontsevich's star product \cite{kontsevich}.
The map $\tilde{\rho}$ of Theorem \ref{thm:cf} has the explicit form
$$
\tilde{\rho}(f) = \sum_\Gamma \tilde{w}_\Gamma D_\Gamma(f).
$$
The sum is over all Kontsevich graphs with one type II\footnote{ Recall that in a Kontsevich graph, there are two kinds of vertices. Type I or ``aerial'' vertices represent one copy of the Poisson structure $\pi$, whereas type II vertices are associated to the functions one intends to multiply.} vertex (associated to $f$). The differential operator $D_\Gamma$ is constructed exactly as it is constructed for Kontsevich's star product. The weights $\tilde{w}_\Gamma$ are given by the integral formula
\[
\tilde{w}_\Gamma = \int_{\tilde{C}_\Gamma} \tilde{\omega}_\Gamma.
\]
Here the $\tilde{C}_\Gamma$, the \emph{Cattaneo-Felder configuration space}, is (a compactification of) the space of all embeddings of the vertex set of $\Gamma$ into the first quadrant, such that the type II vertex is mapped into the real axis. Similar to the Kontsevich case, the weight form $\tilde{\omega}_\Gamma$ is defined as a product of one-forms, one for each edge in the edge set $E(\Gamma)$ of $\Gamma$.
\[
\tilde{\omega}_\Gamma = \pm \bigwedge_{e\in E(\Gamma)} d\phi(z_{e_1}, z_{e_2}).
\]
Here the edge $e$ in the product is understood to point from the vertex $e_1$ that is mapped to $z_{e_1}$ to the vertex $e_2$, that is mapped to $z_{e_2}$.

The precise expression for the angle form $d\phi(z_1, z_2)$ will never be needed, but we will use the following facts about its boundary behaviour in the Appendix:
\begin{enumerate}
\item If $z_1$ lies on the real axis or $z_2$ on the imaginary axis, $d\phi(z_1,z_2)$ vanishes.
\item If $z_1$ and $z_2$ both come close to each other and and a point on the positive real axis, $d\phi(z_1,z_2)$ approaches Kontsevich's angle form.
\item If $z_1$ and $z_2$ both come close to each other and a point on the positive imaginary axis, $d\phi(z_2,z_1)$ approaches Kontsevich's angle form. I.e., $d\phi(z_1,z_2)$ approaches Kontsevich's form after reversal of the edge direction.
\end{enumerate}

\subsection{Construction of the Anomaly}
\label{sec:cfano}
The anomaly $A(f,g)$ in Theorem \ref{thm:cf} can be computed by the formula
$$
A(f,g) = \sum_\Gamma \tilde{w}_\Gamma D_\Gamma (f,g)
$$

Here the sum is over all \emph{anomaly graphs}. Such a graph is a Kontsevich graph, but with a third kind of vertices, which we call \emph{anomalous} or \emph{type III} vertices. An anomaly graph is required to contain at least one such type III vertex. These anomalous vertices have exactly 2 outgoing edges.

The weight $\tilde{w}_\Gamma$ is computed just as the Cattaneo-Felder weight, but with the type III vertices constraint to be mapped to the imaginary axis.

The computation of $D_\Gamma$ also remains the same as before, but one has to specify which bivector fields to associate with the new type III vertices. In local coordinates $x^i$, $i=1,..,n$, the components of this bivector field will be denoted\footnote{The ``$ _a$'' in $\pi_a^{ij}$ is not an index, just a label.}
$$
\epsilon\pi_a^{ij}.
$$
It is in turn given as a sum of graphs.
$$
\pi_a^{ij} = \sum_\Gamma a_{\Gamma} D_\Gamma(x^{i},x^{j})
$$
Here the sum is over all Cattaneo-Felder graphs with 2 type II vertices and $D_\Gamma$ is again defined as in the Cattaneo-Felder case before. However, the weights $a_{\Gamma}$ are computed by the following algorithm:

\begin{enumerate}
\item Delete the type II vertices in $\Gamma$ and all their adjacent edges.
\item Reverse the direction of all edges.
\item Compute the Kontsevich weight of the resulting graph.
\end{enumerate}

\begin{figure}
        \centering
        \includegraphics[width=.8 \columnwidth]{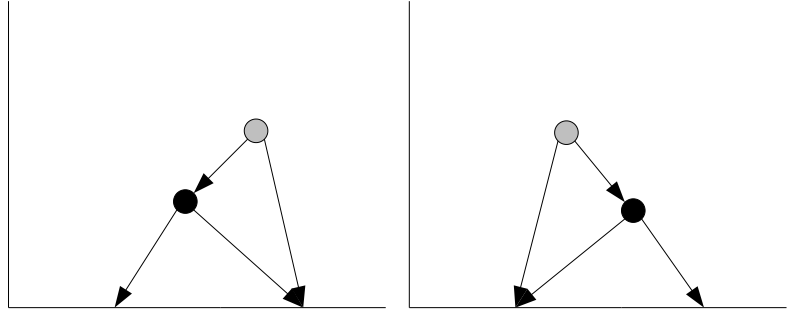}
        \caption{\label{fig:anogra} The two anomaly graphs contributing to the symmetric $\epsilon^3$ part of the anomaly A(f,g). The black vertex is a ``normal'' type I vertex, corresponding to $\pi$. The grey vertex is an anomalous vertex corresponding to $\pi_a$. To this order, $\pi_a=\pi_1$.}
\end{figure}

With this anomalous vertex, one can construct two kinds of graphs that yield $O(\epsilon^3)$-contributions:
\begin{itemize}
\item The graph with only one vertex, which is anomalous. It yields the contribution proportional to $\epsilon \pi_a^{ij}$ to the anomaly.
\item The graphs with one type I and one anomalous vertex as shown in Figure \ref{fig:anogra}. Together, they yield a symmetric contribution to the anomaly.
\end{itemize}

\begin{figure}
        \centering
        \includegraphics[width=.5\columnwidth]{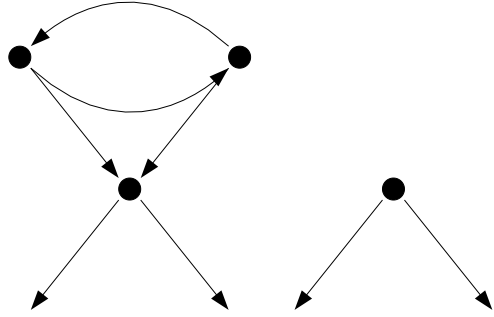} 
        \caption{\label{fig:first_anomally} The graphs accounting for the $\epsilon^2$ and $\epsilon^3$ contributions to the anomaly vertex $\pi_a$. The right graph only contributes if $\pi_\epsilon(0)\neq 0$.} 
\end{figure}

We will use the following notation for the parts of $A(f,g)$ of various orders in $\epsilon$:\footnote{Here and in the following summation over repeated indices is implicit.}
\begin{multline*}
A(f,g) = \epsilon^2 A_2^{ij}(\partial_if)(\partial_jg)+\epsilon^3 A_3^{ij}(\partial_if)(\partial_jg) + \epsilon^3 (\text{symm. in }f,g) + O(\epsilon^4)
\end{multline*}
Here $A_2^{ij}=\pi_1^{ij}$ and $A_3^{ij}$ are antisymmetric. The contribution to $A_3^{ij}$ comes from the two graphs in Figure \ref{fig:first_anomally}, with $\pi's$ attached to the vertices of the left graph and a $\pi_2$ attached to the vertex of the right graph. Note also that if $\pi_\epsilon(0)=0$ to all orders in $\epsilon$, then the contribution of the right graph vanishes and $A_2^{ij}=0$.

\begin{rem}[Linear Poisson structures]
It is easily seen that, if the Poisson structure $\pi_\epsilon=\pi$ is linear, i.e., if $M$ is the dual of a Lie algebra, the anomalous vertex $\pi_a$ vanishes \cite{cattaneo-2004-69}. This is because any contributing graph with $n$ vertices will contain $2n-2$ edges. But for a graph $\Gamma$ with different numbers of vertices and edges, $D_\Gamma=0$ by power counting. Hence a contribution will not arise unless $n=2$. But the weight of the only possible graph with 2 vertices is 0 by one of Kontsevich's lemmatas \cite{kontsevich}, i.e., a reflection argument. 

This also implies that the anomaly can be removed whenever the Poisson structure is linearizable. At least formally, the Poisson structure can be linearized whenever the second Lie algebra cohomology of the Lie algebra defined by the linear order, with values in the symmetric algebra, vanishes. For example, this is true for semisimple Lie algebras. See \cite{weinstein} for details.
\end{rem}

\begin{rem}[Higher order Poisson structures]
As pointed out to the author by A. S. Cattaneo and G. Felder, the anomaly also vanishes whenever the linear order of $\pi$ does. This is shown by power counting: Each vertex in an anomaly graph comes with two edges (derivatives) and is of at least quadratic order. But two edges have to be external and do not contribute derivatives, so the result vanishes when evaluated at $m$.
\end{rem}

\label{sec:cafpoint}

\section*{Acknowledgements}
The author wants to thank his advisor Prof. G. Felder for his continuous support and patient explanation of most of the material presented here. Furthermore, the computation of weights in Appendix \ref{sec:nonvanishing} is mainly due to G. Felder and A. S. Cattaneo, the author merely filled in some details.

\section{Statement of the Main Theorem}

\begin{thm}
\label{thm:main}
There exists a Poisson structure $\pi$ on $\mathbb{R}^n$ vanishing at $0$, s.t. no quantization of $ev_0$ exists for the Kontsevich star product associated to $\pi$.
Furthermore, this also holds for any formal Poisson structure $\pi_\epsilon$ extending $\pi$
$$
\pi_\epsilon=\pi + \epsilon \pi_1 + \cdots
$$
as long as $\pi_\epsilon(0)=0$.
\end{thm}

\begin{rem}
Placing this theorem into a more general context, this means that the $C^\infty(M)$-module structure on a coisotropic submanifold can not always be quantized to a module structure for $(C^\infty(M)[[\epsilon]],\star)$, where $\star$ is the usual Kontsevich product. Hence this theorem answers Main Question 1 negatively. Main Question 2 however, is only partially answered. A complete answer to Main Question 2 we cannot give, only some more hints, see \cite{twdipl}.
\end{rem}

\section{The Proof}
Without loss of generality we can use the following ansatz for $\rho$.
\begin{equation}
\label{equ:rhoansatz}
\rho = \tilde{\rho} + \phi
\end{equation}
Here $\tilde{\rho}$ is the map of Theorem \ref{thm:cf} and $\phi$ has the form
\[
\phi = \epsilon \phi_1 +\epsilon^2 \phi_2 + \cdots.
\]
Our goal, and content of the next sections, is to find the lowest order restrictions on $\phi$ coming from the requirement of $\rho$ being a quantum module map, i.e., eqn. \eref{equ:mod}. Concretely, the requirement is
\begin{equation}
\label{equ:mainmod}
\phi(f\star g)+A(f,g)=\tilde{\rho}(f)\phi(g)+\phi(f)\tilde{\rho}(g)+\phi(f)\phi(g).
\end{equation}
Here we simply inserted \eref{equ:rhoansatz} into \eref{equ:mod} and used \eref{equ:modano}.
To order $\epsilon^0$ this equation is obviously satisfied.

\subsection{Order $\epsilon^1$}
The $\epsilon^1$ part of eqn. \eref{equ:mainmod} reads
\begin{equation}
\label{equ:o1}
\phi_1(fg) = \phi_1(f)g(0) + f(0)\phi_1(g).
\end{equation}
Choosing $f,g$ both constant we see that the zeroth derivative part of $\phi_1$ has to vanish. Choosing $f$ and $g$ both linear the r.h.s. vanishes and hence the second-derivative contribution to $\phi_1$ has to vanish. Picking $f$ quadratic and $g$ linear we see that the third-derivative part of $\rho_1$ must vanish and similarly that all higher derivative parts must vanish as well.

Hence 
\begin{equation}
\label{equ:con1}
\phi_1(f) = D_1^k(\partial_kf)(0)
\end{equation}
for some constants $D_1^k, k=1,..,n$.

\subsection{Order $\epsilon^2$}

We will separately consider the contributions symmetric and antisymmetric in $f$, $g$. 
The antisymmetric contribution reads
\begin{equation}
\label{equ:d1precon}
\phi_1(\aco{f}{g}) = 0.
\end{equation}
Note that if $\pi_\epsilon(0)=0$, then $A_2(f,g)=0$. The l.h.s. of \eref{equ:d1precon} is zero if $f$ or $g$ contains no linear part. Hence it suffices to treat the case where $f$ and $g$ are both linear. Then the equation becomes
\begin{equation}
\label{equ:d1con}
D_1^k(\partial_k\pi^{ij})(0) = 0
\end{equation}
where $\pi^{ij}=\aco{x^i}{x^j}$ are the components of $\pi$ w.r.t. the standard coordinates $\{x^i\}_{i=1,..,n}$.

The symmetric part yields the constraint
\begin{equation}
\label{equ:o2sym}
\phi_2(fg)  =  \phi_1(f)\phi_1(g)+\phi_2(f)g(0)+f(0)\phi_2(g).
\end{equation}

Picking $f,g$ linear we see that the second derivative part of $\phi_2$ must be
$$
\phi_2^{(2)} = \frac{1}{2} D_1^iD_1^j \partial_i\partial_j
$$
Inserting this back into \eref{equ:o2sym} we obtain the same constraint equation for the remaining parts of $\phi_2$ as we had found for $\phi_1$ in eqn. \eref{equ:o1}. By the same logic as there we can hence deduce that
\begin{equation}
\label{equ:con2s}
\phi_2 = \frac{1}{2}D_1^iD_1^j\partial_i\partial_j + D_2^i\partial_i
\end{equation}
for some yet undetermined constants $D_2^i$, $i=1,..,n$.
Here all derivatives are implicitly understood to be evaluated at zero, e.g.
$$
\partial_i\partial_j(f) := (\partial_i\partial_jf)(0).
$$

\begin{rem}
Note that the calculations presented so far are valid for any formal Poisson structure
$$
\pi_\epsilon=\pi + \epsilon \pi_1 + \cdots
$$
as long as it vanishes at 0. I.e., the higher order terms do not contribute to the first two orders in $\epsilon$ of eqn. \eref{equ:mainmod}.
\end{rem}

\subsection{Order $\epsilon^3$}

\begin{figure}
        \centering
        \includegraphics[width=.8\columnwidth]{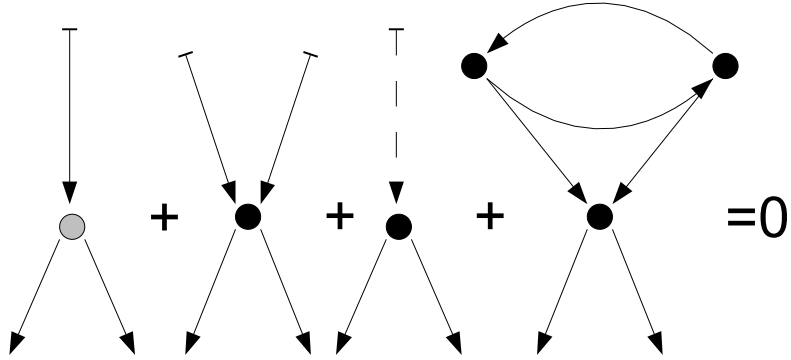}
        \caption{\label{fig:greqn} Eqn. \eref{equ:con3} in terms of graphs. The black vertices symbolize $\pi$, the grey ones $\pi_1$. The solid and dashed arrows without a vertex at their tail symbolize $D_1$ and $D_2$ respectively. The constant coefficients are omitted.} 
\end{figure}

We will only need to consider the antisymmetric part and linear $f=x^i$ and $g=x^j$ in eqn. \eref{equ:mainmod}. In this case the $\epsilon^3$ part of the equation becomes
\begin{equation}
\label{equ:con3}
D_1^k\partial_k\pi_1^{ij} + \frac{1}{2}D_1^kD_1^l\partial_k\partial_l\pi^{ij} + D_2^k\partial_k\pi^{ij}+2A^{ij}_3 = 0.
\end{equation}
The first term is the contribution of the $\epsilon^1$-term in the formal Poisson structure $\pi_\epsilon$. It is absent if we consider $\pi_\epsilon=\pi$.
To derive the above formula we used the following.
\begin{itemize}
\item The r.h.s. of \eref{equ:mainmod} is obviously symmetric in $f$, $g$, hence all contributions come from the l.h.s.
\item The Kontsevich product satisfies
$$
\co{x^i}{x^j}_\star=\epsilon \pi^{ij}+\epsilon^2 \pi_1^{ij} +O(\epsilon^3).
$$
Hence, using \eref{equ:d1con} we obtain
\begin{eqnarray*}
\phi(\co{x^i}{x^j}_\star) &=& \epsilon^3 \phi_1(\pi_1^{ij})+\epsilon^3\phi_2(\pi^{ij}) +O(\epsilon^4)\\
&=& \epsilon^3\left( D_1^k\partial_k\pi_1^{ij} +  \frac{1}{2}D_1^kD_1^l\partial_k\partial_l\pi^{ij} + D_2^k\partial_k\pi^{ij} \right)+O(\epsilon^4)
\end{eqnarray*}
\end{itemize}

The constraint \eref{equ:con3} is displayed graphically in Figure \ref{fig:greqn}.

\subsection{The Counterexample}
\label{sec:thecounterex}
In this section we present a $\pi$ such that there are no constants $D_{1,2}^k$ satisfying \eref{equ:d1con} and \eref{equ:con3} for $\pi_1=0$. This will prove the first part of Theorem \ref{thm:main}.

For this, the following data is needed:
\begin{enumerate}
\item Some finite dimensional semisimple Lie algebra $\alg{g}$ with structure coefficients $g_k^{ij}$ in some basis $\{x^i\}$. We denote by $K^{ij}$ its Killing form and by $K_{ij}$ its inverse. 
\item Some finite dimensional Lie algebra $\alg{h}$ such that its second cohomology group $H^2(\alg{h}, \mathbb{C})\neq \{0\}$. Denote its structure coefficients $h_c^{ab}$.
\item A non-trivial $C\in H^2(\alg{h}, \mathbb{C})$, with coefficients $C^{ab}$ in some basis $\{y^a\}$.
\end{enumerate}

\begin{ex}
The simplest possible choice would be $\alg{g}=\alg{so}(3)$, $\alg{h}=\rn{2}$ as Abelian Lie algebra and $C^{12}=-C^{21}=1$, $C^{11}=C^{22}=0$.
\end{ex}

The Poisson structure will reside in 
$$
\rn{n}\cong (\alg{g}\oplus\alg{h})^* \cong \alg{g}^*\oplus\alg{h}^* 
$$
where $n=\dim\alg{g}+\dim\alg{h}$. We use as coordinate functions the above basis $x^i$ and $y^a$. 
Then one can define
\begin{eqnarray}
\pi^{ij} &=& x^k g^{ij}_k \nonumber \\
\pi^{ai}&=& 0 \label{equ:cexdef} \\
\pi^{ab} &=& y^ch^{ab}_c + \Psi(x)C^{ab} \nonumber
\end{eqnarray}
where $\Psi(x) = K_{ij}x^ix^j$ is the quadratic Casimir in $S^2\alg{g}$.

\begin{lemma}
The above $\pi$ defines a Poisson structure on $\rn{n}$.
\end{lemma}
\begin{proof}
Denote by $\pi^{(1)}$, $\pi^{(2)}$ the linear and quadratic parts of $\pi$ respectively. We need to show that
$$
\co{\pi}{\pi}=\co{\pi^{(1)}}{\pi^{(1)}}+2\co{\pi^{(1)}}{\pi^{(2)}}+\co{\pi^{(2)}}{\pi^{(2)}}=0
$$
where $\co{\cdot}{\cdot}$ denotes the Schouten-Nijenhuis bracket. The linear part of the equation, i.e., $\co{\pi^{(1)}}{\pi^{(1)}}=0$ is satisfied since $\alg{g}\oplus\alg{h}$ is a Lie algebra. The cubic part $\co{\pi^{(2)}}{\pi^{(2)}}=0$ is trivially satisfied since all vector fields $\frac{\partial}{\partial y^a}$ commute with all $x^i$.

The quadratic part $\co{\pi^{(1)}}{\pi^{(2)}}=0$ is equivalent to
$$
\aco{f}{\aco{g}{h}_2}_1+\aco{f}{\aco{g}{h}_1}_2+\textnormal{cycl.} = 0
$$
for all linear $f,g,h\in C^\infty(\rn{n})$. Here $\aco{\cdot}{\cdot}_{1,2}$ are the Poisson brackets of the Poisson structures $\pi^{(1)}$ and $\pi^{(2)}$ respectively. By trilinearity, we can separately consider the following cases.
\begin{itemize}
\item If at least two of the $f,g,h$ are functions of the $x^i$'s only, the expression trivially vanishes since the set of these functions is closed under $\aco{\cdot}{\cdot}_1$, and furthermore $\aco{x^i}{\cdot}_2=0$.

\item If $f=y^a,g=y^b,h=y^c$ the only contributing term is the second, i.e., 
$$
\Psi(x)C^{ad}h^{bc}_d + \textnormal{cycl.}
$$
which vanishes by the cocycle property of $C^{ab}$. 

\item The remaining case is $f=x^i,g=y^a,h=y^b$, leading to (since $\aco{x^i}{y^{a,b}}_{1,2}=0$)
$$
\aco{x^i}{\Psi(x)}_1C^{ab}
$$
which vanishes since the Casimir element is Poisson central.
\end{itemize}
\end{proof}

Knowing that $\pi$ defines a Poisson structure, we can continue the proof of the main theorem. This will be done in two lemmata.

\begin{lemma}
The anomaly $A_3(f,g)$ associated to $\pi$ as in Theorem \ref{thm:cf} is a nonzero multiple of 
$$
C(f,g):= C^{ab}(\partial_af)(0)(\partial_bg)(0).
$$
\end{lemma}
\begin{proof}
The anomaly is given by the left graph, say $\Gamma$, of Figure \ref{fig:first_anomally}.
It will be shown in the Appendix that its weight $\tilde{w}_\Gamma$ is nonzero.
The associated bidifferential operator (applied to functions $f$, $g$) is given by
\begin{eqnarray*}
D_\Gamma(f,g)&=& K^{\alpha\beta}(\partial_\alpha\partial_\beta\pi^{\gamma\delta})(\partial_\gamma f)(\partial_\delta g) \\
&=& K^{ij}(\partial_i\partial_j\Psi(x))C^{ab}(\partial_a f)(\partial_b g) \\
&=& 2K^{ij}K_{ij}C^{ab}(\partial_a f)(\partial_b g) \\
&=& 2(\dim\alg{g}) C^{ab}(\partial_a f)(\partial_b g)
\end{eqnarray*}

Here and in the following we adopt the convention that greek indices refer to a basis
$$
\xi^\alpha=\left\{
\begin{array}{ll}
x^\alpha               & \text{for }\alpha=1,..,\dim\alg{g} \\
y^{\alpha-\dim\alg{g}} & \text{for }\alpha=\dim\alg{g}+1,..,n \\
\end{array} \right.
$$
of $\alg{g}\oplus\alg{h}$, and are summed over $1,..,n$ if repeated. In contrast the roman indices $i$, $j$ label the basis $x^i$ of $\alg{g}$ only and are summed over $1,..,\dim\alg{g}$ if repeated.
Similarly, the roman indices $a$, $b$ refer to the basis $y^a$ of $\alg{h}$ only and are summed over $\dim\alg{g}+1,..,n$.
\end{proof}

\begin{lemma}
For the above $\pi$ and $\pi_1=0$, eqns. \eref{equ:d1con} and \eref{equ:con3} cannot be solved simultaneously.
\end{lemma}
\begin{proof}
Since $\alg{g}$ is semisimple $\co{\alg{g}}{\alg{g}}=\alg{g}$ and eqn. \eref{equ:d1con} implies that 
$$
D_1^i=0,\ i=1,..,\dim(\alg{g}).
$$
But then also 
$$
D_1^\alpha D_1^\beta(\partial_\alpha \partial_\beta\pi^{\gamma\delta})=0
$$
since $\pi$ contains no part quadratic in the $y^a$.
Hence eqn. \eref{equ:con3} becomes
\begin{eqnarray*}
D_2^\alpha (\partial_\alpha \pi^{ij})(0)&=&-2A_3^{ij} = 0 \\
D_2^\alpha (\partial_\alpha \pi^{ab})(0)&=&-2A_3^{ab} \propto C^{ab}
\end{eqnarray*}
Inserting the expression for $\pi$ we obtain from the second equation
$$
D_2^c h^{ab}_c \propto C^{ab}
$$
stating that $C^{ab}$ is a coboundary. Hence, by choice of $C^{ab}$, this equation cannot be solved. Thus the lemma and the first part of Theorem \ref{thm:main} is proven.
\end{proof}

\subsection{A Specialized Counterexample}

We finally turn to the more general case where $\pi_1\neq 0$, but still $\pi_{\epsilon}(0)=0$. The construction in this case runs as above, but we make the special choice $\alg{h}=\alg{k}\oplus\alg{k}$, where $\alg{k}$ is the (unique) non-abelian two dimensional Lie algebra. Its cohomology groups are computed in the Appendix. There is, up to normalization, only one non-trivial cocycle we can pick, namely $\omega$ as defined in eqn. \eref{equ:omegadef} in the Appendix. We will call the resulting Poisson structure $\pi$. The proof of the main Theorem \ref{thm:main} will then be finished by proving the following lemma.

\begin{lemma}
For any formal Poisson structure $\pi_\epsilon$ extending the $\pi$ constructed above, for which $\pi_\epsilon(0)=0$, eqns. \eref{equ:d1con} and \eref{equ:con3} cannot be solved simultaneously.
\end{lemma}
\begin{proof}
By the previous proof it will be sufficient to show that we cannot pick $\pi_1$ and $D_1^a$ such that
\begin{equation}
\label{equ:con3cohom}
D_1^c\partial_c\pi_1^{ab}(0) - \text{(const.)}C^{ab}
\end{equation}
becomes (the coefficients of) an exact element of $H^2(\alg{h},\mathbb{R})$. 
The $\epsilon^1$ part of the condition \eref{equ:piepsjacobi} that $\pi_\epsilon$ is a Poisson structure reads
$$
\co{\pi}{\pi_1} = 0.
$$
Considering only the linear part we have
\begin{equation}
\label{equ:pi1cocycle}
\co{\pi^{(1)}}{\pi_1^{(1)}} = 0.
\end{equation}
where $\pi^{(1)}$, $\pi_1^{(1)}$ are the linear parts of $\pi$ and $\pi_1$ respectively. Note that we used here that the constant part $\pi_1^{(0)}=0$. Eqn. \eref{equ:pi1cocycle} means that $\pi_1^{(1)}$ defines a 2-cocyle of $\alg{g}\oplus\alg{h}$ with values in the adjoint module. 

Equivalently, by projecting on the invariant submodules $\alg{g}$ or $\alg{h}$, one has two 2-cocycles, with values in the $\alg{g}\oplus\alg{h}$-modules $\alg{g}\otimes \mathbb{R}$ and $\mathbb{R}\otimes\alg{h}$ respectively. Here $\mathbb{R}$ is always understood as equipped with the trivial module structure, and $\alg{g}$, $\alg{h}$ with the adjoint structures.

The first 2-cocycle is irrelevant to us since it does not occur in \eref{equ:con3cohom} (since $D_1^i=0$).

The second 2-cocycle defines a cohomology class, say $[\pi_1]$, of
\[
H^2(\alg{g}\oplus\alg{h},\alg{h}).
\]
Eqn. \eref{equ:d1con} means that $D_1^c$ is a cocycle and defines the class
\[
[D_1]\in H^0(\alg{g}\oplus\alg{h}, \alg{h}^*).
\]
The triviality of \eref{equ:con3cohom} implies that their cup product would have to satisfy
\begin{equation}
\label{equ:picoho}
[D_1]\cup [pi_1] \neq 0
\end{equation}

But from the formulas of Knneth and Whitehead it follows that
$$
H^2(\alg{g}\oplus\alg{h},\alg{h})\cong H^0(\alg{g},\mathbb{R})\otimes H^2(\alg{h}, \alg{h}) \cong H^2(\alg{h}, \alg{h}).
$$
Furthermore, as shown in the Appendix, $H^2(\alg{h}, \alg{h})=\{0\}$. Hence $[\pi_1]=0$ and eqn. \eref{equ:picoho} can not be satisified.
Thus the lemma and Theorem \ref{thm:main} is proven.

\end{proof}

\appendix

\section{Cohomology of $\alg{k}$ and $\alg{k}\oplus\alg{k}$}
The Lie algebra $\alg{k}$ is defined as the vector space $\rn{2}$ with the bracket
$$
\co{e_1}{e_2}=e_2
$$
where $e_{1,2}$ are the standard basis vectors.

\begin{lemma}
\label{lem:khom1}
$$
H^0(\alg{k}, \mathbb{R}) \cong H^1(\alg{k}, \mathbb{R}) \cong \mathbb{R} 
$$
All other cohomology groups with values in $\mathbb{R}$ vanish.
A representative of the equivalence class spanning $H^1(\alg{k}, \mathbb{R})$ is 
$$
l: e_1\mapsto 1, \ e_2\mapsto 0.
$$
\end{lemma}
\begin{proof}
It is clear by antisymmetry that $H^n(\alg{k}, \mathbb{C})=\{0\}$ for $n>2$ and also that any 2-cochain is a cocycle.
There is only one 2-cochain (up to a factor) and it is a coboundary since
$$
c: e_1\mapsto 0, \ e_2\mapsto 1
$$
satisfies
$$
c(\co{e_1}{e_2}) = 1.
$$
Finally, any 1-cocycle must vanish on $\co{\alg{k}}{\alg{k}}=\rn{}e_2$. Hence it is clear that the map $l$ defined in the lemma spans the space of 1-cocycles.
\end{proof}

\begin{lemma}
All cohomology groups of $\alg{k}$ with values in $\alg{k}$ vanish, i.e.,
\[
H^p(\alg{k}, \alg{k}) = \{0\} \ \forall p.
\]
\end{lemma}
\begin{proof}
There is no central element in $\alg{k}$, hence $H^0(\alg{k},\alg{k})=\{0\}$.
The cocycle condition for some $l:\alg{k}\rightarrow\alg{k}$ reads
\begin{eqnarray*}
(dl)(e_1,e_2) &=& l(\co{e_1}{e_2})-\co{l(e_1)}{e_2}-\co{e_1}{l(e_2)} \\
&=& l(e_2)-e_2{l_2(e_2)}-e_2l_1(e_1) \\
&=& e_1l_1(e_2)-e_2l_1(e_1) \stackrel{!}{=} 0
\end{eqnarray*}
where $l(\cdot)=e_1l_1(\cdot)+e_2l_2(\cdot)$. Hence $l_1(\cdot)\equiv 0$. But then
\[
l(\cdot) = \co{\cdot}{l_2(e_1)e_2 - l_2(e_2)e_1}
\]
and hence $l$ is exact and $H^1(\alg{k},\alg{k})=\{0\}$. That $H^p(\alg{k},\alg{k})=\{0\}$ for $p\geq 2$ follows as in the proof of the previous lemma.
\end{proof}

We now consider the direct sum $\alg{h}=\alg{k}\oplus\alg{k}$. We denote the standard basis by $e_1$,..,$e_4$. So, e.g., $\co{e_3}{e_4}=e_4$. Knneth's formula and the above lemmata tell us the following:
\begin{itemize}
\item $H^2(\alg{h}, \mathbb{C})$ is spanned by
\begin{equation}
\label{equ:omegadef}
\omega: e_1\wedge e_3\mapsto 1
\end{equation}
with all other components vanishing.
\item $H^2(\alg{h}, \alg{h})=\{0\}$.
\end{itemize}

\section{Nonvanishing of the 2-wheel graph contributing to the anomaly}
\label{sec:nonvanishing}
\begin{figure}
        \centering
                \includegraphics[width=.3\columnwidth]{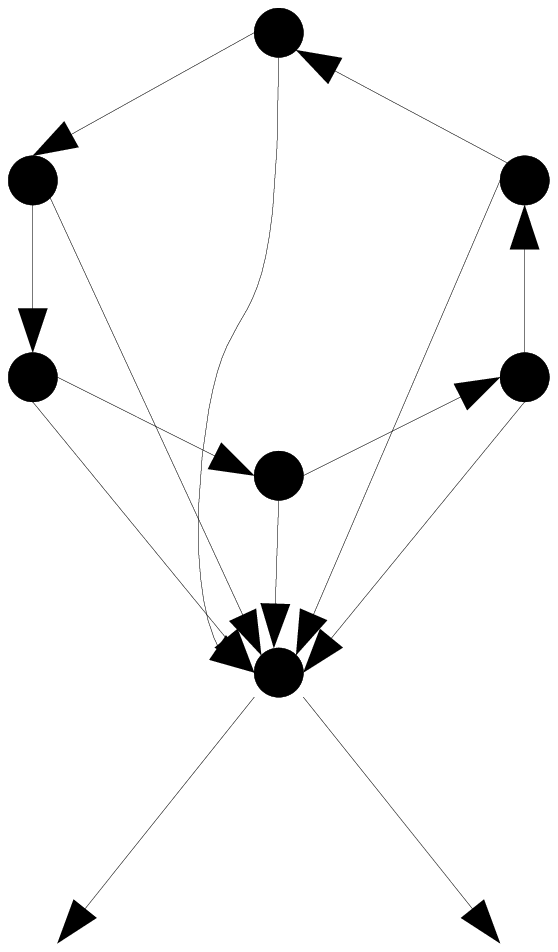}
                \includegraphics[width=.3\columnwidth]{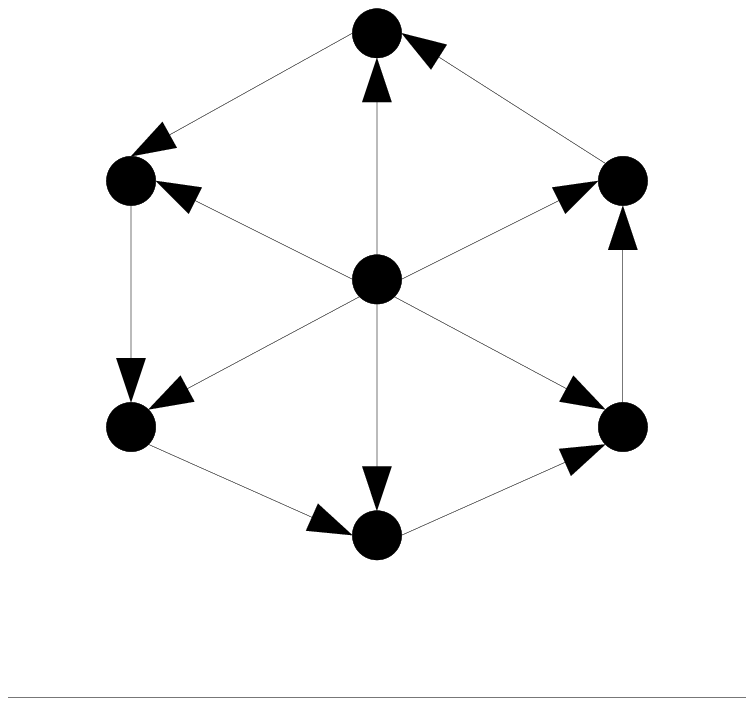}
        \caption{\label{fig:sixwheel} A typical anomaly wheel graph whose weight $c_{2n}$ (here $n=3$) is computed in Theorem \ref{thm:wheelweights}. To get a Cattaneo Felder wheel graph as in the proof of Theorem \ref{thm:wheelweights} one simply removes the two lower edges. The weight is the same as the Kontsevich weight of the graph on the right.}
\end{figure}

One still needs to show that the weight of the left graph in Figure \ref{fig:first_anomally} is nonzero. We will actually compute the weights of all wheel graphs. Instead of defining ``wheel graph'', we refer to Figure \ref{fig:sixwheel}, from which the definition should be clear. To compute the weights, we need the following result interesting in its own right.

\begin{prop}
Let $\alg{g}$ be a Lie algebra and equip its dual space with the canonical Poisson structure. Let $ev_0$ be the evaluation map at zero and $\rho$ its quantization according to Cattaneo and Felder. Let $D:S\alg{g}\rightarrow S\alg{g}$ be the map
$$
D=\left.\det\right.^{\frac{1}{2}}\left(\frac{\sinh(\ad_\partial/2)}{\ad_\partial/2}\right)= \exp\left( \sum_{n\geq 1} \frac{B_{2n}}{4n (2n)!}\tr \left(\ad_\partial^{2n}\right) \right)
$$
with $B_j$ the $j$-th Bernoulli number.\footnote{The map $D$ becomes the Duflo map when composed with $\exp(\tr (ad_\partial)/4)$.}
Then
$$
\rho = ev_0 \circ D^{-1}.
$$  
\end{prop}

\begin{proof}
The map $D$ intertwines the CBH and Kontsevich star products on $S\alg{g}$ (see \cite{kontsevich}, \cite{shoikhet-2000}), i.e.
\begin{equation}
\label{equ:Dintertw}
D (f \star_{CBH} g) = (D f) \star_K (D g) 
\end{equation}
for all $f,g\in S\alg{g}$.

We also have
$$
\rho = ev_0 \circ \exp\left( \sum_{n\geq 1} c_{2n} \tr \left(\ad_\partial^{2n}\right) \right)
$$
for yet undetermined constants $c_{2n}$. Composing both sides of \eref{equ:Dintertw} with $\rho$ and using \eref{equ:modano}\footnote{The anomaly vanishes in this case.} we obtain
$$
(\rho\circ D) (f \star_{CBH} g) = (\rho\circ D)(f) (\rho\circ D)(g).
$$

We want to show recursively that 
$$
d_{2n} := \frac{B_{2n}}{4n (2n)!}+c_{2n}=0
$$
if $d_{2j}=0$ for $j<n$. To do this pick $X\in \alg{g}$ such that $\tr \left(\ad_X^{2n}\right)\neq 0$ and set $f=g=X^n$.
\footnote{One can always find a Lie algebra in which such an $X$ exists. The constants $c_{2n}$ are weights of wheels and do not depend on the Lie algebra.}
Then a straightforward calculation proves the claim.
\end{proof}

\begin{figure}
        \centering
                \includegraphics[width=.4\columnwidth]{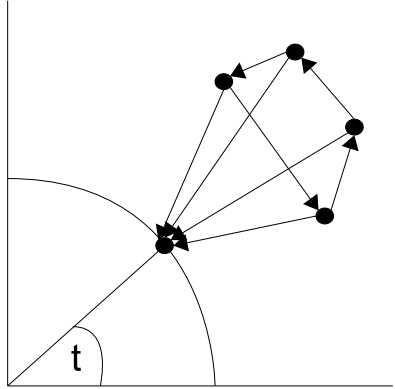}
        \caption{\label{fig:quartcirclemove} A typical configuration in $C_{\Gamma'}$, where $\Gamma$ is the four-wheel graph. Note that the special vertex is confined to the quarter circle to account for the scale invariance.
        }
\end{figure}

\begin{thm}
\label{thm:wheelweights}
The $c_{2n}$ computed in the preceding proof coincide with the weights of the anomaly wheel graphs as depicted in Figure \ref{fig:sixwheel}, up to possibly signs. In particular, the weight of the anomaly graph of Figure \ref{fig:first_anomally} is nonzero.
\end{thm}

\begin{proof}
Pick a Cattaneo Felder wheel graph $\Gamma$. See Figure \ref{fig:sixwheel} for an example. 
Let $\tilde{C}_{\Gamma}$ be the Cattaneo Felder configuration space as in section \ref{sec:cafmod}. 
To divide out the scale invariance we will fix the central vertex of the wheel to lie on the unit quartercircle $\{e^{it};t\in [0,\pi/2]\}$ as depicted in Figure \ref{fig:quartcirclemove}. 

Consider the closed form $\tilde{\omega}_{\Gamma}$ defined on $\tilde{C}_{\Gamma}$ as in section \ref{sec:cafmod}, and compute
\begin{equation}
\label{equ:vanint}
0 = \int_{\tilde{C}_{\Gamma}} d\tilde{\omega}_{\Gamma} = \int_{\partial \tilde{C}_{\Gamma}} \tilde{\omega}_{\Gamma}
\end{equation}
with the help of Stokes' theorem.
There are several boundary strata contributing to the r.h.s. They correspond to center- or non-center-vertices approaching the real axis, imaginary axis, or each other. We divide the strata into the following eight types, each treated separately:\footnote{The readers not familiar with this kind of argument are referred to Kontsevich's proof of his theorem in \cite{kontsevich}.}
\renewcommand{\labelenumi}{(\roman{enumi})}
\begin{enumerate}
\item If all vertices together approach the real axis and each other, the result is 0 by a result of Shoikhet \cite{shoikhet-2000}.
\item If the central vertex approaches the real axis alone, the integral reduces to the integral of $\tilde{\omega}_{\Gamma}$ over $\tilde{C}_{\Gamma}$, yielding the Cattaneo Felder weight $\tilde{c}_\Gamma$.
\item If any subset of vertices approach the real axis and each other, except the two cases before, the result is zero by property 1 in section \ref{sec:cafpoint}. Note that there is always an edge from the collapsing ``cluster'' to the remainder of the graph.
\item If all vertices together approach the imaginary axis and each other the result is the Cattaneo Felder anomaly weight $\tilde{a}_\Gamma$ by property 3 and the algorithm for computing $\tilde{a}_\Gamma$.
\item If any proper subset of vertices appraoch the imaginary axis, the result is zero by property 1.
\item If more than two vertices come close to each other inside the quadrant, the result is zero by a lemma of Kontsevich.
\item If two non-center vertices come close to each other, the result is zero. This is  because both are linked to the center vertex and hence the boundary integrand will contain a wedge product of at least twice the same form, i.e., 0.
\item If any non-center vertex approaches the center vertex, the result is zero by similar reasoning as before. Note that automatically another vertex will be connected twice to the ``cluster'' of the two approaching vertices.  
\end{enumerate}

From this and the vanishing of the integral \eref{equ:vanint} the claim directly follows.

\end{proof}


\end{document}